\documentstyle[11pt,  amscd, amsfonts,leqno]{amsart}

\newtheorem{Theorem}{Theorem}[section]
\newtheorem{Proposition}[Theorem]{Proposition}
\newtheorem{Lemma}[Theorem]{Lemma}

\newtheorem{Example}[Theorem]{Example}

\begin{document}

\title{Extension of plurisubharmonic functions with growth control} 
\author{Dan Coman, Vincent Guedj \and Ahmed Zeriahi}
\thanks{First author is supported by the NSF Grant DMS-0900934}
\date{}
\address{D. Coman: dcoman@@syr.edu, Department of Mathematics, Syracuse
University, Syracuse, NY 13244-1150, USA}
\address{V. Guedj: guedj@@cmi.univ-mrs.fr, Universit\'e Aix-Marseille 1, LATP, 13453 Marseille Cedex 13, FRANCE}
\address{A. Zeriahi: zeriahi@@picard.ups-tlse.fr, Laboratoire Emile
Picard, UMR 5580, Universit\'e Paul Sabatier, 118 route de
Narbonne, 31062 Toulouse Cedex 04, FRANCE}

\pagestyle{myheadings} 

\begin{abstract}
Suppose that $X$ is an analytic subvariety of a Stein manifold $M$ and that $\varphi$ is a plurisubharmonic (psh) function on $X$ which is dominated by a continuous psh exhaustion function $u$ of $M$. Given any number $c>1$, we show that $\varphi$ admits a psh extension to $M$ which is dominated by $cu$ on $M$.

\par We use this result to prove that any $\omega$-psh function on a subvariety of the complex projective space is the restriction of a global $\omega$-psh function, where $\omega$ is the Fubini-Study K\"ahler form.
\end{abstract}

\maketitle

\section*{Introduction}
\par Let $X\subset{\Bbb C}^n$ be a (closed) analytic subvariety. In the case when $X$ is smooth it is well known that a plurisubharmonic (psh) function on $X$ extends to a psh function on ${\Bbb C}^n$ \cite{Sa82} (see also \cite[Theorem 3.2]{BlLe03}). Using different methods, Coltoiu generalized this result to the case when $X$ is singular \cite[Proposition 2]{Col91}.

\smallskip

\par In this article we follow Coltoiu's approach and show that it is possible to obtain extensions with global growth control: 
 
 \medskip
 
\noindent{\bf Theorem A.} 
{\em Let $X$ be an analytic subvariety of a Stein manifold $M$ and let $\varphi$ be a psh function on $X$. Assume that $u$ is a continuous psh exhaustion function on $M$ so that 
$\varphi(z)<u(z)$ for all $z\in X$. Then for every $c>1$ there exists a psh function $\psi=\psi_c$ on $M$ so that $\psi\,|_{_X}=\varphi$ and $\psi(z)<c\max\{u(z),0\}$ for all $z\in M$.} 

\medskip

\par We recall that a function $\varphi:X\to[-\infty,+\infty)$ is called psh if $\varphi\not\equiv-\infty$ on $X$ and if every point $z\in X$ has a neighborhood $U$ in ${\Bbb C}^n$ so that $\varphi=u\,|_{_U}$ for some psh function $u$ on $U$. We refer to \cite{FN} and \cite[section 1]{De85} for a detailed discussion of this notion. We note here that if $\varphi$ is not identically $-\infty$ on an irreducible component $Y$ of $X$ then $\varphi$ is locally integrable on $Y$ with respect to the area measure of $Y$. Let us stress that the more general notion of {\em weakly psh} function is not appropriate for the extension problem (see section \ref{S:alg}).

\medskip

\par We then look at a similar problem on a  compact K\"ahler manifold $V$. Here psh functions have to be replaced by quasiplurisubharmonic (qpsh) ones. Given a K\"ahler form $\omega$, we let
$$
PSH(V,\omega)=\left\{ \varphi \in L^1(V,[-\infty,+\infty)): \, \varphi \text{ upper semicontinuous, } dd^c\varphi \geq -\omega \right\}
$$
denote the set of $\omega$-plurisubharmonic ($\omega$-psh) functions. If $X \subset V$ is an analytic  subvariety, we define similarly the class $PSH(X,\omega\,|_{_X})$ of $\omega$-psh functions on $X$ (see section \ref{S:qpsh} for precise definitions). 

\par By restriction, $\omega$-psh functions on $V$ yield $\omega\,|_{_X}$-psh functions on $X$. Assuming that $\omega$ is a {\it Hodge form}, i.e. a K\"ahler form with integer cohomology class, our second result is that every  $\omega\,|_{_X}$-psh function on $X$ arises in this way.

\medskip

\noindent{\bf Theorem B.}  {\em Let $X$ be a subvariety of a projective manifold $V$ equipped with a Hodge form $\omega$. Then any $\omega\,|_{_X}$-psh function on $X$ is the restriction of an $\omega$-psh function on $V$.}

\medskip

\par Note that in the assumptions of Theorem B there exists a positive holomorphic line bundle $L$ on $V$ whose first Chern class $c_1(L)$ is represented by $\omega$. In this case the $\omega$-psh functions are in one-to-one correspondence with the set of (singular) positive metrics of $L$ (see \cite{GZ}). Thus an alternate formulation of Theorem B is the following:  
 
\medskip

\noindent{\bf Theorem B'.}  {\em Let $X$ be a subvariety of a projective manifold $V$ and $L$ be an ample line bundle on $V$. Then any (singular) positive metric of $L\,|_{_X}$ is the restriction of a (singular) positive metric of $L$ on $V$.}

\medskip
 
\par Recall that it is possible to regularize qpsh functions on ${\Bbb P}^n$, since it is a homogeneous manifold. Hence Theorem B has the following immediate corollary:

\medskip

\noindent{\bf Corollary C.}  {\em Let $X$ be a subvariety of a projective manifold $V$ equipped with a Hodge form $\omega$. If $\varphi \in PSH(X,\omega\,|_{_X})$ then there exists a sequence of smooth functions $\varphi_j\in PSH(V,\omega)$ which decrease pointwise on $V$ so that $\lim\,\varphi_j=\varphi$ on $X$.}

\medskip

\par When $X$ is smooth this regularization result is well known to hold even when the cohomology class of $\omega$ is not integral (see \cite{De92}, \cite{BK07}). 

\par Corollary C allows to show that the singular K\"ahler-Einstein currents constructed in \cite{EGZ} have {\em continuous} potentials, a result that has been obtained recently in \cite{EGZ2} by completely different methods (see also \cite{DZ} for partial results in this direction).

\medskip

\par We prove Theorem A in section \ref{S:pfTA}. The compact setting is considered in section \ref{S:qpsh}, where Theorem B is derived from Theorem A. In section \ref{S:alg} we discuss the special situation when $X$ is an algebraic subvariety of ${\Bbb C}^n$. As an application of Theorem B, we give a characterization of those psh functions in the Lelong class ${\mathcal L}(X)$ which admit an extension in the Lelong class ${\mathcal L}({\Bbb C}^n)$ (see section \ref{S:alg} for the necessary definitions). In particular, we give simple examples of algebraic curves $X\subset{\Bbb C}^2$ and of functions $\eta\in{\mathcal L}(X)$ which do not have extensions in ${\mathcal L}({\Bbb C}^2)$.

\section{Proof of Theorem A}\label{S:pfTA}

\par The following proposition will allow us to reduce the proof of Theorem A to the case $M={\Bbb C}^n$. 
We include its short proof for the convenience of the reader. 

\begin{Proposition}\label{P:Sad} Let $V$ be a complex submanifold of ${\Bbb C}^N$ and $u$ be a continuous psh exhaustion function on $V$. Then there exists a continuous psh exhaustion function $\widetilde u$ on ${\Bbb C}^N$ so that $\widetilde u\,|_{_V}=u$.
\end{Proposition}

\begin{pf} The argument is very similar to the one of Sadullaev (\cite{Sa82},\cite[Theorem 3.2]{BlLe03}). By \cite{Siu76}, there exists an open neighborhood $W$ of $V$ in ${\Bbb C}^N$ and a holomorphic retraction $r:W\to V$.
We can find an open neighborhood $U$ of $V$ so that $U\subset W$ and $\|r(z)-z\|<2$ for every $z\in U$. Indeed, if $B(p,r)$ denotes the open ball in ${\Bbb C}^N$ centered at $p$ and of radius $r$, then $U_p=r^{-1}(B(p,1))\cap B(p,1)$ is an open neighborhood of $p\in V$, and we let $U=\bigcup_{p\in V}U_p$. Since $u$ is a continuous psh exhaustion function on $V$, it follows that the function $u(r(z))$ is continuous psh on $U$
and $\lim_{z\in U,\|z\|\to+\infty}u(r(z))=+\infty$. 

\par It is well known that there exist entire functions $f_0,\dots,f_N$, so that $V=\{z\in{\Bbb C}^N:\,f_k(z)=0,\;0\leq k\leq N\}$ (see \cite[p.63]{Ch}). The function $\rho=\log(\sum|f_k|^2)$ is psh on ${\Bbb C}^N$ and $V=\{\rho=-\infty\}$. 

\par Let $D$ be an open set so that $V\subset D\subset\overline D\subset U$. Since $\rho$ is continuous on ${\Bbb C}^N\setminus V$, we can find a convex increasing function $\chi$ on $[0,+\infty)$ which verifies for every $R\geq0$ the following two properties:

\par $(i)$ $\chi(R)>R-\rho(z)$ for all $z\in{\Bbb C}^N\setminus D$ with $\|z\|=R$.

\par $(ii)$ $\chi(R)>u(r(z))-\rho(z)$ for all $z\in\partial D$ with $\|z\|=R$. \\
Then  
$$\widetilde u(z)=\left\{\begin{array}{ll}\max\{u(r(z)),\chi(\|z\|)+\rho(z)\},\;{\rm if}\;z\in D,\\
\chi(\|z\|)+\rho(z),\;{\rm if}\;z\in{\Bbb C}^N\setminus D,
\end{array}
\right.$$
is a continuous psh exhaustion function on ${\Bbb C}^N$ and $\widetilde u=u$ on $V$.
\end{pf}

\par Employing the methods of Coltoiu \cite{Col91} we now construct psh extensions with growth control  over bounded sets in ${\Bbb C}^n$.

\begin{Proposition}\label{P:Col} Let $\chi$ be a psh function on a subvariety $X\subset{\Bbb C}^n$ and let $v$ be a continuous psh function on ${\Bbb C}^n$ with $\chi<v$ on $X$. If $R>0$, there exists a psh function $\widetilde\chi=\widetilde\chi_R$ on ${\Bbb C}^n$ so that $\widetilde\chi\,|_{_X}=\chi$ and 
$\widetilde\chi(z)<v(z)$ for all $z\in{\Bbb C}^n$ with $\|z\|\leq R$.
\end{Proposition}

\begin{pf} We use a similar argument to the one in the proof of Proposition 2 in \cite{Col91}. Consider the subvariety $A=(X\times{\Bbb C})\cup({\Bbb C}^n\times\{0\})\subset{\Bbb C}^{n+1}$, and let 
$$D=\{(z,w)\in X\times{\Bbb C}:\,\log|w|+\chi(z)<0\}\cup({\Bbb C}^n\times\{0\})\subset A.$$ 
Since $D\cap(X\times{\Bbb C})$ is Runge in $X\times{\Bbb C}$, it follows that $D$ is Runge in $A$. Let 
$$K=\{(z,w)\in{\Bbb C}^{n+1}:\,\rho(z,w)=\max\{\log^+(\|z\|/R),\log|w|+v(z)\}\leq0\}.$$
Since $v$ is continuous, $\rho$ is a continuous psh exhaustion function on ${\Bbb C}^{n+1}$, so $K$ is a polynomially convex compact set. As $\chi<v$ on $X$, we have $K\cap A\subset D$. By \cite[Theorem 3]{Col91} there exists a Runge domain $\widetilde D\subset{\Bbb C}^{n+1}$, with $\widetilde D\cap A=D$ and $K\subset\widetilde D$. Let $\delta(z,w)$ denote the distance from $(z,w)\in\widetilde D$ to 
$\partial\widetilde D$ in the $w$-direction. Since $\widetilde D$ is pseudoconvex, $-\log\delta$ is psh on $\widetilde D$ (see e.g. \cite[Proposition 9.2]{FS}). Hence $\widetilde\chi(z)=-\log\delta(z,0)$ is psh on ${\Bbb C}^n$, as ${\Bbb C}^n\times\{0\}\subset\widetilde D$. Since $\widetilde D\cap A=D$, it follows that $\widetilde\chi\,|_{_X}=\chi$. Moreover, $K\subset\widetilde D$ implies that $\widetilde\chi(z)<v(z)$ for all $z\in{\Bbb C}^n$ with $\|z\|\leq R$.
\end{pf}

\par The proof of Theorem A proceeds like this. Given a partition $${\Bbb C}^n=\bigcup\{m_{j-1}<u\leq m_j\},$$ where $m_j\nearrow+\infty$, we apply Proposition \ref{P:Col} inductively to construct an extension dominated in each ``annulus" $\{m_{j-1}<u\leq m_j\}$ by $\gamma_ju$, where $\gamma_j>1$ is an increasing sequence defined in terms of the $m_j$'s. Theorem A will follow by showing that it is possible to choose $\{m_j\}$ rapidly increasing so that $\lim \gamma_j$ is arbitrarily close to 1. 

\par We fix next an increasing sequence $\{m_j\}_{j\geq-1}$ so that 
$$m_{-1}=m_0=0<m_1<m_2<\dots,\:\{u<m_1\}\neq\emptyset,\;m_j\nearrow+\infty.$$
Define inductively a sequence $\{\gamma_j\}_{j\geq0}$, as follows:
\begin{equation}\label{e:gamma}
\gamma_0=1,\;\;\gamma_j(m_j-m_{j-1})=\gamma_{j-1}(m_j-m_{j-2})+1\;{\rm for}\;j\geq1.
\end{equation}
Clearly, $\gamma_j>\gamma_{j-1}>1$ for all $j>1$. 

\begin{Proposition}\label{P:mj} 
Let $X,\,\varphi,\,u$ be as in Theorem A with $M={\Bbb C}^n$, and let $\{m_j\}$, $\{\gamma_j\}$ be as above. There exists a psh function $\psi$ on ${\Bbb C}^n$ so that $\psi\,|_{_X}=\varphi$ and for all $z\in{\Bbb C}^n$ we have 
$$\psi(z)<\left\{\begin{array}{ll}\gamma_j u(z),\;{\rm if}\;m_{j-1}<u(z)\leq m_j,\;j\geq2,\\
\gamma_1\max\{u(z),0\},\;{\rm if}\;u(z)\leq m_1.
\end{array}\right.$$
\end{Proposition}

\begin{pf} We introduce the sets 
$$D_j=\{z\in{\Bbb C}^n:\,u(z)<m_j\}\;,\;\;K_j=\{z\in{\Bbb C}^n:\,u(z)\leq m_j\}.$$
Since $u$ is a continuous psh exhaustion function, $K_j$ is a compact set. Let
$$\rho_j=\gamma_j\max\{u-m_{j-1},0\}-j,\;j\geq0.$$
Then $\rho_j$ is psh on ${\Bbb C}^n$ and (\ref{e:gamma}) implies that
\begin{equation}\label{e:rho1}
\rho_j(z)=\rho_{j-1}(z)\;{\rm if}\;u(z)=m_j,\;j\geq1.
\end{equation}
We claim that 
\begin{equation}\label{e:rho2}
\rho_j(z)\geq u(z)\;{\rm if}\;z\in{\Bbb C}^n\setminus D_j,\;j\geq0.
\end{equation}
Indeed, since $\gamma_j\geq1$ and using (\ref{e:gamma}) we obtain
\begin{eqnarray*}
\rho_j(z)-u(z)&=&(\gamma_j-1)u(z)-\gamma_jm_{j-1}-j\geq(\gamma_j-1)m_j-\gamma_jm_{j-1}-j\\
&=&(\gamma_{j-1}-1)m_j-\gamma_{j-1}m_{j-2}-j+1\\
&\geq&(\gamma_{j-1}-1)m_{j-1}-\gamma_{j-1}m_{j-2}-(j-1).
\end{eqnarray*}
So $x_j:=(\gamma_j-1)m_j-\gamma_jm_{j-1}-j\geq x_0=0$, and (\ref{e:rho2}) is proved.

\par Let $\varphi_j=\max\{\varphi,-j\}$. We construct by induction on $j\geq1$ a sequence of continuous psh functions $\psi_j$ on ${\Bbb C}^n$ with the following properties:
\begin{eqnarray}
\label{e:i}&& \psi_j(z)>\varphi_j(z)\;{\rm for}\;z\in X\;,\;\;\int_{X\cap K_{j-1}}(\psi_j-\varphi_j)<2^{-j}.\\
\label{e:ii}&& \psi_j(z)\geq\rho_j(z)\;{\rm for}\;z\in D_j\;,\;\; 
\psi_j(z)=\rho_j(z)\;{\rm for}\;z\in{\Bbb C}^n\setminus D_j.\\
\label{e:iii}&& \psi_j(z)<\psi_{j-1}(z)\;{\rm for}\;z\in K_{j-1},\;{\rm where}\;\psi_0=\rho_0=\max\{u,0\}.
\end{eqnarray}
Here the integral in (\ref{e:i}) is with respect to the area measure on each irreducible component, i.e.
$$\int_{X\cap K}f:=\sum\int_{Y\cap K}f\,\beta^{\dim Y},$$
where the sum is over all irreducible components $Y$ of $X$ which intersect $K$ and $\beta$ is the standard 
K\"ahler form on ${\Bbb C}^n$. (Note that this is a finite sum.)

\par Assume that the function $\psi_{j-1}$ is constructed with the desired properties. We construct $\psi_j$ by applying Proposition \ref{P:Col} with $\chi=\varphi_j$ and $v=\psi_{j-1}$. (If $j=1$, $\psi_1$ is constructed in the same way by applying Proposition \ref{P:Col} with $\chi=\varphi_1$ and $v=\psi_0$.) By (\ref{e:i}), 
$\varphi_j\leq\varphi_{j-1}<\psi_{j-1}$ on $X$ (and for $j=1$, clearly $\varphi_1<\psi_0$ on $X$). Therefore Proposition \ref{P:Col} yields a psh function $\widetilde\varphi_j$ on ${\Bbb C}^n$ so that 
$\widetilde\varphi_j\,|_{_X}=\varphi_j$ and $\widetilde\varphi_j<\psi_{j-1}$ on $K_j$. Using the standard regularization of $\widetilde\varphi_j$ and the dominated convergence theorem (as $\varphi_j\geq-j$) we obtain a continuous psh function $\widetilde\psi_j$ on ${\Bbb C}^n$ which verifies
$$\widetilde\psi_j(z)>\varphi_j(z)\;{\rm for}\;z\in X\;,\;\;\int_{X\cap K_j}(\widetilde\psi_j-\varphi_j)<2^{-j}.$$
Moreover, since $\psi_{j-1}$ is continuous, we can ensure by the Hartogs lemma that we also have  
$\widetilde\psi_j(z)<\psi_{j-1}(z)$ for $z\in K_j$.

\par We now define 
$$\psi_j(z)=\left\{\begin{array}{ll}\max\{\widetilde\psi_j(z),\rho_j(z)\},\;{\rm if}\;z\in D_j,\\
\rho_j(z),\;{\rm if}\;z\in{\Bbb C}^n\setminus D_j.\end{array}\right.$$
By (\ref{e:ii}) and (\ref{e:rho1}) we have $\widetilde\psi_j<\psi_{j-1}=\rho_{j-1}=\rho_j$ on $\partial D_j$ (for $j=1$, recall that $\psi_0=\rho_0$ by definition). So $\psi_j$ is a continuous psh function on ${\Bbb C}^n$ which verifies (\ref{e:ii}). On $X\setminus D_j$ we have by (\ref{e:rho2}) that $\psi_j=\rho_j\geq u>\varphi_j$, while on $X\cap D_j$, $\psi_j\geq\widetilde\psi_j>\varphi_j$. Since $\rho_j=-j\leq\varphi_j<\widetilde\psi_j$ on $X\cap K_{j-1}$, we see that $\psi_j=\widetilde\psi_j$ on $X\cap K_{j-1}$ so
$$\int_{X\cap K_{j-1}}(\psi_j-\varphi_j)\leq\int_{X\cap K_j}(\widetilde\psi_j-\varphi_j)<2^{-j}.$$
Hence $\psi_j$ verifies (\ref{e:i}). Finally, we have by (\ref{e:ii}), $\rho_j=-j<\rho_{j-1}\leq\psi_{j-1}$ on $K_{j-1}$ (and for $j=1$, $\rho_1=-1<\psi_0=0$ on $K_0$). Since $\widetilde\psi_j<\psi_{j-1}$ on $K_j$ we conclude that $\psi_j<\psi_{j-1}$ on $K_{j-1}$, so (\ref{e:iii}) is verified. 

\par So we have constructed a sequence of continuous psh functions $\psi_j$ on ${\Bbb C}^n$ verifying properties (\ref{e:i})-(\ref{e:iii}). Since $\bigcup_{j\geq1} D_j={\Bbb C}^n$, we have by (\ref{e:iii}) that the function $$\psi(z)=\lim_{j\to\infty}\psi_j(z)$$ is well defined and psh on ${\Bbb C}^n$. As $\ldots<\psi_{j+2}<\psi_{j+1}<\psi_j$ on $K_j$, it follows that $\psi<\psi_j$ on $K_j$. 

\par Suppose now that $z\in K_j\setminus D_{j-1}$, for some $j\geq2$, so $m_{j-1}\leq u(z)\leq m_j$. By the above construction and property (\ref{e:ii}), we have 
$$\widetilde\psi_j(z)<\psi_{j-1}(z)=\rho_{j-1}(z)\Longrightarrow\psi(z)<\psi_j(z)\leq\max\{\rho_{j-1}(z),\rho_j(z)\}\leq\gamma_ju(z).$$
Similarly, for $z\in K_1$ we have 
$$\psi(z)<\psi_1(z)\leq\max\{\rho_0(z),\rho_1(z)\}\leq\gamma_1\max\{u(z),0\}.$$
Hence $\psi$ satisfies the desired global upper estimates on ${\Bbb C}^n$.

\par Property (\ref{e:i}) implies that $\psi(z)\geq\varphi(z)$ for every $z\in X$. Let $K$ be a compact in ${\Bbb C}^n$ and $Y$ be an irreducible component of $X$ so that $\varphi\,|_{_Y}\not\equiv-\infty$. By (\ref{e:i}) we have that for all $j$ sufficiently large
$$0\leq\int_{Y\cap K}(\psi_j-\varphi)=\int_{Y\cap K}(\psi_j-\varphi_j)+\int_{Y\cap K}(\varphi_j-\varphi)\leq
2^{-j}+\int_{Y\cap K}(\varphi_j-\varphi).$$
Hence by dominated convergence, $\int_{Y\cap K}(\psi-\varphi)=0$, which shows that $\psi=\varphi$ on $Y$.

\par Assume now that $Y$ is an irreducible component of $X$ so that $\varphi\,|_{_Y}\equiv-\infty$. Then using (\ref{e:i}) and the monotone convergence theorem we conclude that 
$$\int_{Y\cap K}\psi=\lim_{j\to\infty}\int_{Y\cap K}\psi_j=\lim_{j\to\infty}\left(\int_{Y\cap K}(\psi_j-\varphi_j)+\int_{Y\cap K}\varphi_j\right)=-\infty,$$
so $\psi\,|_{_Y}\equiv-\infty$. Therefore $\psi=\varphi$ on $X$, and the proof is finished. \end{pf}

\medskip

\noindent{\em Proof of  Theorem A.} We consider first the case $M={\Bbb C}^n$. Fix $c>1$. We define inductively a sequence $\{m_j\}$ with the following properties: $m_{-1}=m_0=0<m_1$, $\{u<m_1\}\neq\emptyset$, and for $j\geq1$, $m_j>m_{j-1}$ is chosen large enough so that 
$$a_j=\frac{m_{j-1}-m_{j-2}+1}{m_j-m_{j-1}}\leq\frac{\log c}{2^j}\;.$$
Since $\gamma_j\geq\gamma_0=1$ we have by (\ref{e:gamma}),
$$\gamma_j(m_j-m_{j-1})\leq\gamma_{j-1}(m_j-m_{j-2}+1)\Longrightarrow
\gamma_j\leq\gamma_{j-1}(1+a_j).$$
Thus
$$\gamma_j<\gamma=\prod_{j=1}^\infty(1+a_j)\;,\;\;\log\gamma\leq\sum_{j=1}^\infty a_j\leq\log c.$$

\par Let $\psi=\psi_c$ be the psh extension of $\varphi$ provided by Proposition \ref{P:mj} for this sequence $\{m_j\}$. Then for every $z\in{\Bbb C}^n$ we have
$$\psi(z)<\gamma\max\{u(z),0\}\leq c\max\{u(z),0\}.$$

\medskip

\par Assume now that $M$ is a Stein manifold of dimension $n$. Then $M$ can be properly embedded in ${\Bbb C}^{2n+1}$, hence we may assume that $M$ is a complex submanifold of ${\Bbb C}^{2n+1}$ (see e.g. \cite[Theorem 5.3.9]{Hor}). Proposition \ref{P:Sad} implies the existence of a continuous psh exhaustion function $\widetilde u$ on ${\Bbb C}^{2n+1}$ so that $\widetilde u=u$ on $M$. By what we already proved, given $c>1$  there exists a psh function $\widetilde\psi$ on ${\Bbb C}^{2n+1}$ which extends $\varphi$ and such that 
$\widetilde\psi<c\max\{\widetilde u,0\}$ on ${\Bbb C}^{2n+1}$. We let $\psi=\widetilde\psi\,|_{_M}$. $\Box$

\medskip

\par We end this section by noting that some hypothesis on the growth of $u$ is necessary in Theorem A. Indeed, suppose that $X$ is a submanifold of ${\Bbb C}^n$ for which there exists a non-constant negative psh function $\varphi$ on $X$. Then any psh extension of $\varphi$ to ${\Bbb C}^n$ cannot be bounded above. However, by Theorem A, given any $\varepsilon>0$ there exists a psh function $\psi=\psi_\varepsilon$ so that $\psi\,|_{_X}=\varphi$ and $\psi(z)<\varepsilon\log^+\|z\|$ on ${\Bbb C}^n$.

\section{Extension of qpsh functions}\label{S:qpsh}
\par Let $V$ be a compact K\"ahler manifold equipped with a K\"ahler form $\omega$.
We let $PSH(V,\omega)$ denote the set of $\omega$-psh functions on $V$.
These are upper semicontinuous functions $\varphi \in L^1(V,[-\infty,+\infty))$
such that $\omega+dd^c \varphi \geq0$, where $d=\partial+\overline\partial$ and
$d^c=\frac{1}{2 \pi i}(\partial-\overline\partial)$.
We refer the reader to \cite{GZ} for basic properties of $\omega$-psh functions.

\par Let $X$ be an analytic subvariety of $V$. Recall that an upper semicontinuous function $\varphi:X\to[-\infty,+\infty)$ is called $\omega\,|_{_X}$-psh if $\varphi\not\equiv-\infty$ on $X$ and if there exist an open cover $\{U_i\}_{i\in I}$ of $X$ and psh functions $\varphi_i,\rho_i$ defined on $U_i$, where $\rho_i$ is smooth and $dd^c\rho_i=\omega$, so that $\rho_i+\varphi=\varphi_i$ holds on $X\cap U_i$, for every $i\in I$. Moreover, $\varphi$ is called {\em strictly} $\omega\,|_{_X}$-psh if it is $(1-\varepsilon)\omega\,|_{_X}$-psh for some small $\varepsilon>0$. The current $\omega\,|_{_X}+dd^c\varphi$ is then called a K\"ahler current on $X$ (see \cite[section 5.2]{EGZ}). We denote by $PSH(X,\omega\,|_{_X})$, resp. $PSH^+(X,\omega\,|_{_X})$, the class of $\omega\,|_{_X}$-psh, resp. strictly $\omega\,|_{_X}$-psh functions on $X$.

\par Every $\omega$-psh function $\varphi$ on $V$ yields, by restriction, an $\omega\,|_{_X}$-psh function $\varphi\,|_{_X}$ on $X$, as soon as $\varphi\,|_{_X} \not\equiv -\infty$. The question we address here is whether 
this restriction operator is surjective. In other words, is there equality
$$PSH(X,\omega\,|_{_X}) \stackrel{?}{=}PSH(V,\omega)\,|_{_X} .$$

\subsection{The smooth case}
We start with the elementary observation that smooth strictly $\omega$-psh functions can easily be  extended.

\begin{Proposition}
Let $V$ be a compact K\"ahler manifold equipped with a K\"ahler form $\omega$,
and let $X$ be a complex submanifold of $V$. Then
$$
PSH^+(X,\omega\,|_{_X}) \cap {\mathcal C}^{\infty}(X,{\Bbb R})=
\left(PSH^+(V,\omega) \cap {\mathcal C}^{\infty}(V,\Bbb R)  \right)\,|_{_X}.
$$
\end{Proposition}

\par We include a proof for the convenience of the reader, although this is probably part of the ``folklore" (see e.g. \cite{Sch} for the case where $\omega$ is a Hodge form).

\begin{pf}
Let $\varphi \in {\mathcal C}^{\infty}(X,\Bbb R)$ be such that 
$(1-\varepsilon)\omega\,|_{_X}+dd^c \varphi \geq0$ on $X$, for some $\varepsilon>0$. We first choose $\tilde{\varphi}$ to be any smooth extension of $\varphi$ to $V$. Consider
$$
\psi:=\tilde{\varphi}+A \chi\,\text{dist}(\cdot,X)^2,
$$
where $\chi$ is a test function supported in a small neighborhood of $X$
and such that $\chi \equiv 1$ near $X$. Here ${dist}$ is any Riemannian distance on $V$, for instance
the distance associated to the K\"ahler metric $\omega$.
Then $\psi$ is yet another smooth extension of $\varphi$ to $V$,
which now satisfies $(1-\varepsilon/2)\omega+dd^c \psi \geq0$ near $X$,
if $A$ is chosen large enough.

\par The function $\log(\text{dist}(\cdot,X)^2)$ is well defined and qpsh in a neighborhood of $X$. Let $\chi$ be a test function supported in this neighborhood so that $\chi\equiv1$ near $X$. The function $u=\chi\log(\text{dist}(\cdot,X)^2)$ is $N\omega$-psh on $V$ for a large integer $N$. Moreover, $\exp(u)$ is smooth and $X=\{u=-\infty\}$. Replacing $\omega$ by $N\omega$, $\varphi$ by $N\varphi$, and $\psi$ by $N\psi$, we may assume that $N=1$. Set now
$$
\psi_C:=\frac{1}{2}\,\log \left[ e^{2\psi}+e^{u+C} \right].
$$
This again is a smooth extension of $\varphi$, and a straightforward computation
yields 
$$
dd^c \psi_{C} \geq  \frac{2e^{2\psi}dd^c \psi+e^{u+C} dd^c u}{2(e^{2\psi}+e^{u+C})}\;.
$$
Hence 
$$
\left(1-\frac{\varepsilon}{2}\right)\omega+dd^c \psi_C \geq 
\frac{2e^{2\psi}\left[\left(1-\frac{\varepsilon}{2}\right)\omega+dd^c \psi\right]+(1-\varepsilon)e^{u+C} \omega}{2(e^{2\psi}+e^{u+C})}\geq0,
$$
if $C$ is chosen large enough.
\end{pf}

\par This proof breaks down when $\varphi$ is singular and hence a different approach is needed. We consider in the next section the particular case when $\omega$ is  a Hodge form.

\subsection{Proof of Theorem B}
We assume here that $\omega$ is a {\it Hodge form}, i.e. that the cohomology class $\{\omega\}$ belongs to
$H^2(V,\Bbb Z)$ (more precisely to the image of $H^2(V,\Bbb Z)$ in $H^2(V,\Bbb R)$ under the mapping induced by the inclusion $\Bbb Z \hookrightarrow\Bbb R$). We prove the following more precise version of Theorem B.

\begin{Theorem}\label{T:extcomp}
Let $X$ be a subvariety of a projective manifold $V$ equipped with a Hodge form $\omega$. If $\varphi \in PSH(X,\omega\,|_{_X})$ then given any constant $a>0$ there exists $\psi\in PSH(V,\omega)$ so that $\psi\,|_{_X}=\varphi$ and  $\max_V\psi<\max_X\varphi+a$.
\end{Theorem}

\par In the assumptions of Theorem \ref{T:extcomp} there exists a positive holomorphic line bundle $L$ on $V$ whose first Chern class $c_1(L)$ is represented by $\omega$. By Kodaira's embedding theorem $L$ is ample,  hence for large $k$ there exists an embedding $\pi:V\hookrightarrow{\Bbb P}^n$ such that $L^k=\pi^* {\mathcal O}(1)$.

\par Replacing $\omega$ by $k \omega$, $\varphi$ by $k \varphi$, we can assume that $L={\mathcal O}(1)$, $V$ is an algebraic submanifold of the complex projective space ${\Bbb P}^n$, and $\omega=\omega_{FS}\,|_{_V}$ is the Fubini-Study K\"ahler form. Hence $X$ is an algebraic subvariety of ${\Bbb P}^n$, and Theorem \ref{T:extcomp} follows if we show that $\omega_{FS}$-psh functions on $X$ extend to $\omega_{FS}$-psh functions on ${\Bbb P}^n$.

\smallskip
 
\par Therefore we assume in the sequel that $X\subset V={\Bbb P}^n$ and $\omega$ is the Fubini-Study K\"ahler form on ${\Bbb P}^n$. Let $[z_0:\ldots:z_n]$ denote the homogeneous coordinates. Without loss of generality, we may assume that they are chosen so that no coordinate hyperplane $\{z_j=0\}$ contains any irreducible component of $X$. 

\par Let 
$$\theta(z)=\log\frac{\max\{|z_0|,\dots,|z_n|\}}{\sqrt{|z_0|^2+\ldots+|z_n|^2}}\,,\;z=[z_0:\ldots:z_n]\in{\Bbb P}^n.$$
This is an $\omega$-psh function and for all $z\in{\Bbb P}^n$,
$$-m\leq\theta(z)\leq0\,,\;\;{\rm where}\;m=\log\sqrt{n+1}.$$

\par We start by noting that Theorem A yields special subextensions of $\omega$-psh functions on $X$. 

\begin{Lemma}\label{L:subext}
Let $\varepsilon\geq0$ and $u$ be a continuous $(1+\varepsilon)\omega$-psh function on ${\Bbb P}^n$ so that $u(z)\leq0$ for all $z\in{\Bbb P}^n$. If $c>1$ and $\varphi$ is an $\omega$-psh function on $X$ so that $\varphi<u$, then there exists a $c\omega$-psh function $\psi$ on ${\Bbb P}^n$ so that 
$$
\frac{1}{c}\,\psi(z)\leq\frac{1}{1+\varepsilon}\,u(z),\;\forall z\in{\Bbb P}^n,
$$
and 
$$
\psi(z)=\varphi(z)+(c-1)\theta(z)+(c-1)\min_{\zeta\in{\Bbb P}^n}u(\zeta),\;\forall z\in X.
$$
\end{Lemma}

\begin{pf} Let 
$$M=-\min_{\zeta\in{\Bbb P}^n}u(\zeta)\geq0.$$
We work first in an affine chart $\{z_j=1\}\equiv{\Bbb C}^n$. Let $X_j=X\cap\{z_j=1\}$ and let $\rho_j\geq0$ be the potential of $\omega$  in this chart with $\rho_j(0)=0$. Then $\varphi+\rho_j$ is psh on $X_j$ and since $u\leq0$, 
$$\varphi+\rho_j+M<u+\rho_j+M\leq\frac{1}{1+\varepsilon}\,u+\rho_j+M\;{\rm on}\;X_j.$$ 
Note that $(1+\varepsilon)^{-1}u+\rho_j+M\geq0$ is a continuous psh exhaustion function on ${\Bbb C}^n$. Theorem A yields a psh function $\widetilde\psi$ on ${\Bbb C}^n$ so that 
$$\widetilde\psi<\frac{c}{1+\varepsilon}\,u+c\rho_j+cM\;{\rm on}\;{\Bbb C}^n\;,\;\;\widetilde\psi=\varphi+\rho_j+M\;{\rm on}\;X_j.$$

\par The function $\psi_j=\widetilde\psi-c\rho_j-cM$ extends uniquely to a $c\omega$-psh function on ${\Bbb P}^n$ which verifies $$\psi_j\leq\frac{c}{1+\varepsilon}\,u\;\;{\rm on}\;{\Bbb P}^n.$$ 
Moreover on $X\cap\{z_j=1\}$ we have 
$$\psi_j=\varphi-(c-1)\rho_j-(c-1)M=\varphi+(c-1)\theta_j-(c-1)M,$$
where
$$\theta_j(z)=\log\frac{|z_j|}{\sqrt{|z_0|^2+\ldots+|z_n|^2}}\;.$$
Hence $\psi_j=-\infty$ on $X\cap\{z_j=0\}$. 

\par We finally let $\psi=\max\{\psi_0,\ldots,\psi_n\}$. This is a $c\omega$-psh function on ${\Bbb P}^n$ which verifies the desired conclusions, since $\theta=\max\{\theta_0,\ldots,\theta_n\}$.
\end{pf}

\medskip

\par\noindent{\em Proof of Theorem \ref{T:extcomp}}. Fix $a>0$. Replacing $\varphi$ by $\varphi-\max_X\varphi-a$ we may assume that $\max_X\varphi=-a$. We will show that there exists a sequence of smooth $\omega$-psh functions $\varphi_j$ on ${\Bbb P}^n$ which decrease pointwise on ${\Bbb P}^n$ to a negative $\omega$-psh function $\psi$ so that $\psi=\varphi$ on $X$.

\smallskip

\par Let $X'$ be the union of the irreducible components $W$ of $X$ so that $\varphi\,|_{_W}\not\equiv-\infty$. We first construct by induction on $j\geq1$ a sequence of numbers $\varepsilon_j\searrow0$ and a sequence of negative smooth $(1+\varepsilon_j)\omega$-psh functions $\psi_j$ on ${\Bbb P}^n$ so that for all $j\geq2$
$$\frac{\psi_j}{1+\varepsilon_j}<\frac{\psi_{j-1}}{1+\varepsilon_{j-1}}\;\;{\rm on}\;{\Bbb P}^n\;,\;\;\psi_{j-1}>\varphi\;{\rm on}\;X\;,\;\;\int_{X'}(\psi_j-\varphi)<\frac{1}{j}\;,\;\;
\int_{W}\psi_j<-j\,,$$
for every irreducible component $W$ of $X$ where $\varphi\,|_{_W}\equiv-\infty$. Here the integrals are with respect to the area measure on each irreducible component $X_j$ of $X$, i.e.
$$\int_Xf:=\sum_{X_j}\int_{X_j}f\,\omega^{\dim X_j}.$$

\par Let $\varepsilon_1=1$, $\psi_1=0$, and assume that $\varepsilon_{j-1},\,\psi_{j-1}$, where 
$j\geq2$, are constructed with the above properties. Since $\varphi<\psi_{j-1}\,|_{_X}$ and the latter is continuous on the compact set $X$, we can find $\delta>0$ so that $\varphi<\psi_{j-1}-\delta$ on $X$. 

\par Let $c>1$. By Lemma \ref{L:subext}, there exists a $c\omega$-psh function $\psi_c$ so that 
$$\frac{\psi_c}{c}\leq\frac{\psi_{j-1}-\delta}{1+\varepsilon_{j-1}}\;\;{\rm on}\;{\Bbb P}^n\;,\;\;\psi_c=\varphi+(c-1)\theta-(c-1)M_{j-1}\;\;{\rm on}\;X,$$
where
$$M_{j-1}=\delta-\min_{\zeta\in{\Bbb P}^n}\psi_{j-1}(\zeta)\geq0.$$
We can regularize $\psi_c$ on ${\Bbb P}^n$: there exists a sequence of smooth $c\omega$-psh functions decreasing to $\psi_c$ on ${\Bbb P}^n$. Therefore we can find a smooth $c\omega$-psh function $\psi'_c$ on ${\Bbb P}^n$ so that 
$$\frac{\psi'_c}{c}<\frac{\psi_{j-1}-\frac{\delta}{2}}{1+\varepsilon_{j-1}}\;{\rm on}\;{\Bbb P}^n,\;\;\psi'_c>\varphi+(c-1)\theta-(c-1)M_{j-1}\geq\varphi-(c-1)(m+M_{j-1})\;\;{\rm on}\;X.$$
By dominated, resp. monotone convergence, we can in addition ensure that 
$$\int_{X'}(\psi'_c-\varphi)\leq\int_{X'}(\psi'_c-\varphi-(c-1)\theta+(c-1)M_{j-1})<c-1,$$
$$\int_{W}\psi'_c<-j-(c-1)(m+M_{j-1})|W|,$$
for every irreducible component $W$ of $X$ where $\varphi\,|_{_W}\equiv-\infty$. Here $|W|$ denotes the (projective) area of $W$.

\par Now let $\psi''_c=\psi'_c+(c-1)(m+M_{j-1})$. Then on ${\Bbb P}^n$ we have 
$$\frac{\psi''_c}{c}<\frac{\psi_{j-1}-\frac{\delta}{2}}{1+\varepsilon_{j-1}}+\frac{(c-1)(m+M_{j-1})}{c}<\frac{\psi_{j-1}}{1+\varepsilon_{j-1}}-\frac{\delta}{4}+(c-1)(m+M_{j-1}).$$
Moreover, $\psi''_c>\varphi$ on $X$ and 
\begin{eqnarray*}
\int_{X'}(\psi''_c-\varphi)&=&\int_{X'}(\psi'_c-\varphi)+(c-1)(m+M_{j-1})|X'|\\
&<&(c-1)(1+m|X'|+M_{j-1}|X'|)\;,\\
\int_{W}\psi''_c&=&\int_{W}\psi'_c+(c-1)(m+M_{j-1})|W|<-j\;,
\end{eqnarray*}
for every irreducible component $W$ of $X$ where $\varphi\,|_{_W}\equiv-\infty$.

\par We take $c=1+\varepsilon_j$ and $\psi_j=\psi''_c$, where $\varepsilon_j>0$ is so that 
$$\varepsilon_j<\varepsilon_{j-1}/2\;,\;\;\varepsilon_j(m+M_{j-1})<\frac{\delta}{4}\;,\;\;\varepsilon_j(1+m|X'|+M_{j-1}|X'|)<\frac{1}{j}\;.$$ 
Then $\varepsilon_j,\,\psi_j$ have the desired properties. 

\smallskip

\par We conclude that  $\varphi_j=(1+\varepsilon_j)^{-1}\psi_j$ is a decreasing sequence of smooth negative $\omega$-psh function on ${\Bbb P}^n$, so that $\varphi_j>(1+\varepsilon_j)^{-1}\varphi>\varphi$ on $X$. Hence
$\psi=\lim_{j\to\infty}\varphi_j$ is a negative $\omega$-psh function on ${\Bbb P}^n$ and $\psi\geq\varphi$ on $X$. Note that 
$$\int_{X'}(\varphi_j-\varphi)=\frac{1}{1+\varepsilon_j}\int_{X'}(\psi_j-\varphi)-\frac{\varepsilon_j}{1+\varepsilon_j}\int_{X'}\varphi<\frac{1}{j}-\frac{\varepsilon_j}{1+\varepsilon_j}\int_{X'}\varphi\;,$$
$$\int_{W}\varphi_j=\frac{1}{1+\varepsilon_j}\int_{W}\psi_j<-\frac{j}{2}\;,$$
for every irreducible component $W$ of $X$ where $\varphi\,|_{_W}\equiv-\infty$. It follows that $\psi=\varphi$ on $X$ and the proof of Theorem \ref{T:extcomp} is finished. $\Box$

\section{Algebraic subvarieties of ${\Bbb C}^n$}\label{S:alg}
\par If $X$ is an analytic subvariety of ${\Bbb C}^n$ and $\gamma$ is a positive number, we denote by 
${\mathcal L}_\gamma(X)$ the {\em Lelong class} of psh functions $\varphi$ on $X$ which verify 
$\varphi(z)\leq\gamma\log^+\|z\|+C$ for all $z\in X$, where $C$ is a constant that depends on $\varphi$. We let ${\mathcal L}(X)={\mathcal L}_1(X)$. By Theorem A, functions $\varphi\in{\mathcal L}(X)$ admit a psh extension in each class ${\mathcal L}_\gamma({\Bbb C}^n)$, for every $\gamma>1$. \footnote{If $X$ is algebraic this result is claimed in \cite[Proposition 3.3]{BlLe03}, but there is a gap in their proof.}

\par We  assume in the sequel that $X$ is an {\em algebraic} subvariety of ${\Bbb C}^n$ and address the question whether it is necessary to allow the arbitrarily small additional growth. More precisely, is it true that
$$
{\mathcal L}(X) \stackrel{?}{=} {\mathcal L}({\Bbb C}^n)\,|_{_X},
$$
i.e. is every psh function with logarithmic growth on $X$ the restriction of a
globally defined psh function with logarithmic growth? We will give a criterion for this to hold, but show that in general this is not the case.

\subsection{Extension preserving the Lelong class}\label{SS:Lel}
Consider the standard embedding  
$$z\in{\Bbb C}^n\hookrightarrow[1:z]\in{\Bbb P}^n,$$ 
where  $[t:z]$ denote the homogeneous coordinates on ${\Bbb P}^n$. Let $\omega$ be the Fubini-Study K\"ahler form and let 
$$\rho(t,z)=\log\sqrt{|t|^2+\|z\|^2}$$
be its logarithmically homogeneous potential on ${\Bbb C}^{n+1}$. 

\par We denote by $\overline X$ the closure of $X$ in ${\Bbb P}^n$, so $\overline X$ is an algebraic subvariety of ${\Bbb P}^n$. It is well known that the class $PSH({\Bbb P}^n,\omega)$ is in one-to-one
correspondence with the Lelong class ${\mathcal L}({\Bbb C}^n)$ (see \cite{GZ}). Let us look at the connection between $\omega$-psh functions on $\overline X$ and the class ${\mathcal L}(X)$. 

\par The mapping 
$$F_X:PSH(\overline X,\omega\,|_{_{\overline X}})\longmapsto{\mathcal L}(X),\;(F_X\varphi)(z)=\rho(1,z)+\varphi([1:z]),$$
is well defined and injective. However, it is in general not surjective, as shown by Examples \ref{E:nec1} and \ref{E:nec2} that follow.

\par Conversely, a function $\eta\in{\mathcal L}(X)$ induces an upper semicontinuous function $\widetilde\eta$ on $\overline X$ defined in the obvious way:
$$\widetilde\eta([t:z])=\left\{\begin{array}{ll}
\eta(z)-\rho(1,z),\;\hspace{24mm}{\rm if}\;t=1,\;z\in X, \\ \\ 
\displaystyle\limsup_{[1:\zeta]\to[0:z],\zeta\in X}(\eta(\zeta)-\rho(1,\zeta)),\;{\rm if }\;t=0,\;[0:z]\in\overline X\setminus X.
\end{array}\right.$$
The function $\widetilde\eta$ is in general only {\em weakly $\omega$-psh on $\overline X$}, i.e. it is bounded above on $\overline X$ and it is $\omega\,|_{_{\overline X_r}}$-psh on the set $\overline X_r$ of regular points of $\overline X$. This notion is in direct analogy to that of {\em weakly psh} function on an analytic variety (see \cite[section 1]{De85}). We do not pursue it any further here. 

\par Note that $\eta\in F_X\left(PSH(\overline X,\omega\,|_{_{\overline X}})\right)$ if and only if $\widetilde\eta\in PSH(\overline X,\omega\,|_{_{\overline X}})$. The following simple characterization is a consequence of Theorem B. 

\begin{Proposition}\label{P:char} Let $\eta\in{\mathcal L}(X)$. The following are equivalent:

\par (i) There exists $\psi\in{\mathcal L}({\Bbb C}^n)$ so that $\psi=\eta$ on $X$.

\par (ii) $\widetilde\eta\in PSH(\overline X,\omega\,|_{_{\overline X}})$.

\par (iii) For every point $a\in\overline X\setminus X$ the following holds: if $(X_j,a)$ are the irreducible components of the germ $(\overline X,a)$ then the value 
$$\limsup_{X_j\ni[1:\zeta]\to a}(\eta(\zeta)-\rho(1,\zeta))$$ 
is independent of $j$. 

\smallskip

\par In particular,  if the germs  $(\overline X,a)$ are irreducible for all points $a\in\overline X\setminus X$ then ${\mathcal L}(X) ={\mathcal L}({\Bbb C}^n)\,|_{_X}$.
\end{Proposition}

\begin{pf} Assume that $(i)$ holds. It follows that $\widetilde\eta=\varphi\,|_{_{\overline X}}$, where  
$$\varphi([t:z]):=\left\{\begin{array}{ll}
\psi(z)-\rho(1,z),\;\hspace{29mm}{\rm if}\;t=1,\\
\limsup_{[1:\zeta]\to[0:z]}(\psi(\zeta)-\rho(1,\zeta)),\;{\rm if }\;t=0,
\end{array}\right.$$
is an $\omega$-psh function on ${\Bbb P}^n$. Hence $\widetilde\eta\in PSH(\overline X,\omega\,|_{_{\overline X}})$.

\par Conversely, if $(ii)$ holds then by Theorem B there exists an $\omega$-psh function $\varphi$ on ${\Bbb P}^n$ which extends $\widetilde\eta$. Hence $\psi(z)=\rho(1,z)+\varphi([1:z])$ is an extension of $\eta$ and $\psi\in{\mathcal L}({\Bbb C}^n)$.

\par The equivalence of $(ii)$ and $(iii)$ follows easily from \cite[Theorem 1.10]{De85}.
\end{pf}

\subsection{Explicit examples}\label{SS:ex}
In view of section \ref{SS:Lel}, it is easy to construct examples of algebraic curves $X\subset{\Bbb C}^2$ and functions in ${\mathcal L}(X)$ which do not admit an extension in ${\mathcal L}({\Bbb C}^2)$. We write $z=(x,y)\in{\Bbb C}^2$.

\begin{Example}\label{E:nec1} Let $X=\{y=0\}\cup\{y=1\}\subset{\Bbb C}^2$ and $\eta\in{\mathcal L}(X)$, where
$$\eta(z)=\left\{\begin{array}{ll}
\rho(1,z),\;\;\;\;\;\;\;{\rm if}\;z=(x,0),\\
\rho(1,z)+1,\;{\rm if}\;z=(x,1).
\end{array}\right.$$
The function $\widetilde\eta$ is not $\omega$-psh on $\overline X=\{y=0\}\cup\{y=t\}$, hence $\eta$ does not have an extension in ${\mathcal L}({\Bbb C}^2)$. Indeed, the maximum principle is violated along $\{y=0\}$ near the point $a=[0:1:0]$, since $\widetilde\eta([t:1:0])=0$ for $t\neq0$, while $\widetilde\eta([t:1:t])=1$.
\end{Example}

\par With a little more effort we can give an example as above where $X$ is an irreducible curve. Let ${\Bbb C}^\star={\Bbb C}\setminus\{0\}$. 

\begin{Example}\label{E:nec2} Let $X\subset{\Bbb C}^2$ be the irreducible cubic with equation $xy=x^3+1$. Then $$\overline X=\{[t:x:y]\in{\Bbb P}^2:\,xyt=x^3+t^3\},\;\overline X=X\cup\{a\},\;a=[0:0:1].$$ 
The germ $(\overline X,a)$ has two irreducible components $X_1,\,X_2$, both are smooth at $a$, $X_1$ being tangent to the line $\{x=0\}$, and $X_2$ to the line $\{t=0\}$.

\par Note that in fact $X\subset{\Bbb C}^\star\times{\Bbb C}$ is the graph of the rational function $y=x^2+x^{-1}$, $x\in{\Bbb C}^\star$. If $(x,y)\in X$ and $x\to 0$ then $(x,y)\to a$ along $X_1$, while as $x\to\infty$ then $(x,y)\to a$ along $X_2$. The function 
$$u(x,y)=\max\{-\log|x|,2\log|x|+1\}$$ 
is psh in ${\Bbb C}^\star\times{\Bbb C}$. It is easy to check that $\eta:=u\,|_{_X}\in{\mathcal L}(X)$ and 
$$\limsup_{X_1\ni[1:\zeta]\to a}(\eta(\zeta)-\rho(1,\zeta))=0\;,\;\;\limsup_{X_2\ni[1:\zeta]\to a}(\eta(\zeta)-\rho(1,\zeta))=1.$$
Hence $\eta$ does not admit an extension in ${\mathcal L}({\Bbb C}^2)$.
\end{Example}

\par We conclude this section with an example of a cubic $X$ in ${\Bbb C}^2$ and a psh function on $X$ of the form  $\eta=\log|P|$, where $P$ is a polynomial, so that $\eta$ admits a ``transcendental" extension with exactly the same growth, but small additional growth is necessary if we look for an ``algebraic" extension.

\begin{Proposition}
Let $X=\{x=y^3\}$ and $\eta(x,y)=\log|1+y|$, so $\eta\,|_{_X}\in{\mathcal L}_{1/3}(X)$.

\par Given $k\geq1$, there is a polynomial $Q_k(x,y)$ of degree $k+1$ so that $Q_k(y^3,y)=(y+1)^{3k}$. In particular,
$\psi_k=\frac{1}{3k} \log|Q_k|\in {\mathcal L}_{(k+1)/3k}({\Bbb C}^2)$
is an extension of $\eta\,|_{_X}$.

\par There exists no polynomial $Q(x,y)$ of degree $k$ so that 
$Q(y^3,y)=(y+1)^{3k}$. However, $\eta\,|_{_X}$ has an extension in ${\mathcal L}_{1/3}({\Bbb C}^2)$.
\end{Proposition}

\begin{pf} We construct $Q_k$ by replacing $y^3$ by $x$ in the polynomial 
$$(y+1)^{3k}=\sum_{j=0}^{3k}{3k\choose j}y^j.$$
Since $j=3[j/3]+r_j$, $r_j\in\{0,1,2\}$, it follows that 
$$
Q_k(x,y)=\sum_{j=0}^{3k}{3k\choose j}x^{[j/3]}y^{r_j}=3kx^{k-1}y^2+l.d.t.\;.
$$
 
\par We now check that there is no polynomial $Q(x,y)$ of degree $k$ so that 
$Q(y^3,y)=(y+1)^{3k}$. Indeed, if $Q(x,y)=\sum_{j+l\leq k}c_{jl}x^jy^l$ then
$$Q(y^3,y)=c_{k0}y^{3k}+c_{k-1,1}y^{3k-2}+l.d.t.$$
does not contain the monomial $y^{3k-1}$. 

\par Note that $\overline X=\{xt^2=y^3\}=X\cup\{a\}$, where $a=[0:1:0]$, so the germ $(\overline X,a)$ is irreducible. Proposition \ref{P:char} implies that $\eta\,|_{_X}$ has an extension in ${\mathcal L}_{1/3}({\Bbb C}^2)$.
\end{pf}

\par We conclude with some remarks regarding our last example. If $X$ is an algebraic subvariety of ${\Bbb C}^n$ and $f$ is a holomorphic function on $X$, $f$ is said to have polynomial growth if there is an integer $N(f)$ and a constant $A$ so that 
$$|f(z)|\leq A(1+\|z\|)^{N(f)},\;\;\forall\,z\in X.$$
Then it is well known that there exists a polynomial $P$ of degree at most $N(f)+\varepsilon(X)$ so that $P\,|_{_X}=f$, where $\varepsilon(X)>0$ is a constant depending only on $X$ (see e.g. \cite{Bj74} and references therein). However, if $\overline X\subset{\Bbb P}^N$ is irreducible at each of its points at infinity then by Proposition \ref{P:char} the psh function $\eta=N(f)^{-1}\log|f|\in{\mathcal L}(X)$ has a psh extension in the Lelong class ${\mathcal L}({\Bbb C}^n)$. 

\par On the other hand, Demailly \cite{De79} has shown that in the case of the transcendental curve $X=\{e^x+e^y=1\}$ any holomorphic function $f$ on $X$, of polynomial growth, has a polynomial extension of the same degree to ${\Bbb C}^n$. Hence it is natural to ask if for this curve one has that ${\mathcal L}(X)={\mathcal L}({\Bbb C}^n)\,|_{_X}$.

 \end{document}